\documentclass[12pt,a4paper]{amsart}
\usepackage[T1]{fontenc}

\usepackage{amssymb}
\usepackage{color}
\usepackage{graphicx,color}
\usepackage{amsmath, amssymb, graphics}
\usepackage[latin1]{inputenc}

\theoremstyle{theorem}

\newtheorem*{classification}{The classification}
\newtheorem*{method}{The method}

\theoremstyle{definition}
\newtheorem*{definition}{Definition}
\newtheorem*{remark}{Remark}

\begin{document}

\title{How does a topologist classify the letters of the alphabet?}
\markright{Classify the letters of the alphabet}
\author{Rafael L\'opez}
 \address{Departamento de Geometr\'{\i}a y Topolog\'{\i}a\\
Universidad de Granada\\
18071 Granada, Spain\\}
\thanks{Partially supported by MEC-FEDER
 grant no. MTM2011-22547 and Junta de Andaluc\'{\i}a grant no. P09-FQM-5088.}
 \email{ rcamino@ugr.es}

\maketitle

\begin{figure}[hbtp]
\begin{center}\includegraphics{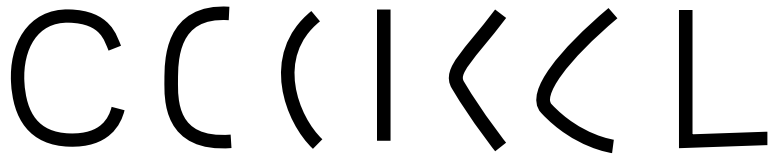}\end{center}
\caption{The transformations between the letters {\sf C}, {\sf I} and {\sf L} by stretching and bending  show that all are homeomorphic.}\label{fig1}
\end{figure}

The letters of the alphabet, when written on a piece of paper, may be viewed as subsets of the plane.  Imagine that  each  letter is made by an elastic material in such a way that you can deform it by hand, by  stretching, contracting or twisting, but you can never tear, fold, glue or cut.  The question that we propose you is whether it is possible to deform a letter in another one by the operations listed above. On the contrary, if you can not do it, how to find a method to decide that there does not exist such a deformation. This type of problems are topological.  A topologist thinks an object made by rubber and he studies those properties that remain unchanged under transformations by stretching and contracting. In our  context of Euclidean plane ${\mathbb R}^2$,    if $X,Y\subset{\mathbb R}^2$ are two subsets, we say that $X$ is homeomorphic to $Y$, and we write $X\cong Y$, if there exists a one-to-one mapping $\phi:X\rightarrow Y$  which is bi-continuous. This transformation is called a {\it  homeomorphism} and for a topologist, $X$  and $Y$ are identical. For example, the letters {\sf C}, {\sf I} and {\sf L} are homeomorphic such as it is illustrated in Fig. \ref{fig1}.

This problem on the classification of the alphabet appeared as an exercise proposed by Gustave Choquet in \cite[p. 21]{ch}, where he intends ``{\it to convey the intuitive content of homeomorphness}". Even Wikipedia gathers the problem as ``an introductory exercise'' \cite{wiki}. In fact, in the figure that Choquet shows in \cite{ch}, the letters  {\sf C I J L M N S U V W Z} are all homeomorphic.  However, letters as {\sf A} and {\sf E} seem not be equivalent. The work of distinguishing the  letters of the alphabet is not discussed in \cite{ch} and we need the use of the concept of  \emph{topological invariant} to assert that, indeed,   {\sf A} and {\sf E} are not homeomorphic.

Two remarks. First, we view the letter written on paper as a curve in the plane, that is, a $1$-dimensional objets  and thus, our letters have zero thickness. Second, each person has its own   alphabet and one  writes a letter with a special style. An enthusiastic reader is invited to work with her/his own alphabet. In typography,  there are many fonts of a typeface, with a specific  style, italicization and design.  In order to unify the variety of possibilities, we here  consider the sans  serif font of \TeX, which it is obtained  with   the \TeX \ command \verb|\sf|. Furthermore, we will only classify uppercase  letters such as it is depicted in Fig. \ref{fig2}.

In what follows we shall write a letter of the alphabet in sans serif font if it is viewed as a subset of ${\mathbb R}^2$ and if the letter denotes a mathematical symbol, then we write it in roman font. For example,  {\sf X} indicates the letter of the alphabet and $X$ a  symbol.

\begin{figure}[hbtp]
\begin{center}{\Large{\sf A B C D E F G H I J K L M N O P Q R S T U V W X Y Z}}\end{center}
\caption{The uppercase  letters of the alphabet written  in  sans serif font by \TeX.}\label{fig2}
\end{figure}
  In order to distinguish  two letters of the alphabet, we will use a topological property that satisfies a letter but not the other one. Roughly speaking, the idea is the following.     Let $X, Y\subset{\mathbb R}^2$ be two letters as subsets of Euclidean plane. Take a point $p\in X$ and we remove $p$ from $X$, that is, we consider the subset $X-\{p\}\subset {\mathbb R}^2$. When we delete $p$ from $X$, it is possible that $X$ splits into different pieces just as when a mirror falls and breaks.  Denote by $O(p)$ the number of pieces $X-\{p\}$. If $X\cong Y$, then there must be a point $q\in Y$ with the same   property, that is,  $O(q)=O(p)$. Furthermore, if for each natural number $n\in{\mathbb N}$, $N(n,X)$ is the number of points of $X$ with $O(p)=n$, then $N(n,X)=N(n,Y)$ for all $n\in{\mathbb N}$. The work will consist in the next steps: (i) define, in a simple way, a piece  of a subset of ${\mathbb R}^2$; (ii) give a  method to count the number of pieces of a set; (iii) compute $N(n,X)$ for all letters of the alphabet and (iv) use the computations of (iii) to get the classification of the letters.

\section{Path connectedness in Euclidean space.}

We need some basics and terminology from topology  \cite{mu}. Denote by ${\mathbb R}$ the set of real number and by ${\mathbb R}^m={\mathbb R}\times\stackrel{\stackrel{m}{\smile}}{\ldots}\times  {\mathbb R}$ the $m$-space. We introduce the Euclidean topology of ${\mathbb R}^m$  thanks to the Euclidean distance. If $x=(x_1,\ldots,x_n)$, $y=(y_1,\ldots,y_n)$ are two points of ${\mathbb R}^m$, the Euclidean distance between $x$ and $y$ is
$$d(x,y)=\sqrt{(x_1-y_1)^2+\ldots+(x_n-y_n)^2}.$$
If $x\in{\mathbb R}^m$ and $r>0$, the $m$-ball centered at $x$ of radius $r$ is the set
$B(x;r)=\{y\in{\mathbb R}^m:d(y,x)<r\}$. For example, an $1$-ball is an open interval of ${\mathbb R}$ and a $2$-ball is a round disc of ${\mathbb R}^2$. A subset $O\subset{\mathbb R}^m$ is an\emph{ open set} of the Euclidean topology if for all $x\in O$ there is a $m$-ball centered at $x$ and contained in $O$. By considering all $x\in O$ and if $r_x>0$ is the corresponding radius  such that $B(x;r_x)\subset O$, then $\cup_{x\in O} B(x;r_x)\subset O$. Hence   we deduce $O=\cup_{x\in O}B(x;r_x)$ and this proves that an open  of ${\mathbb R}^m$ is the union of $m$-balls.  This topology is carried on  a subset $X \subset {\mathbb R}^m$ by saying that a subset $A\subset X$ is called \emph{open in $X$} if it is the intersection of an open set $O$ of ${\mathbb R}^m$ with $X$, that is, $A=X\cap O$.

The topological notion that we need for the classification of the alphabet is  path connectedness.
A {\it path} in ${\mathbb R}^m$ is a continuous map $\alpha:[0,1]\rightarrow {\mathbb R}^m$. If $p=\alpha(0)$ and $q=\alpha(1)$, we say that $p$ and $q$ are the initial and final points of $\alpha$, respectively, and that $\alpha$ joins $p$ with $q$.  Intuitively a path in Euclidean space may view as a (continuous) trajectory of a point. A set $X\subset {\mathbb R}^m$ is said {\it path connected} if for any $p,q\in X$ there is a path {\it in $X$} joining $p$ and $q$. The path connected property is topological in such a way that  is any set homeomorphic to  a path connected set is path connected. An easy property  is that the union of a collection of  path connected sets  $\{X_i:i\in I\}$ with non-trivial intersection  is path connected because if $p_0\in \cap_{i\in I}X_i$, then  two points of $\cup_{i\in I}X_i$ can be continuously linked by a path throughout the point $p_0$. Examples of path connected spaces are: an interval of ${\mathbb R}$, a straight-line of ${\mathbb R}^m$, a circle, a sphere and the very space ${\mathbb R}^m$.

\begin{figure}[hbtp]
\begin{center}
\includegraphics{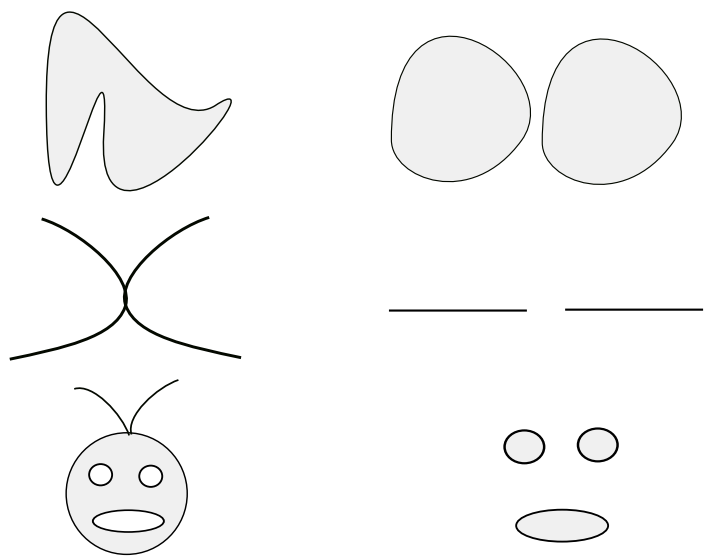}
\end{center}
\caption{On the left column, three path connected sets; on the right column, three non-path connected sets}\label{fig3}
\end{figure}

If $X\subset {\mathbb R}^m$ is not path connectedness, given $p\in X$, we call the {\it path component} of $p$, and denoted by $C_p$,  the largest path connected subset of $X$ containing $p$. Then the component $C_p$ is the set of all points joining with $p$  by using paths of $X$ and this establishes a partition of $X$. Each component is then a piece of $X$ and the number of path components is topological, so, if $\phi:X\rightarrow Y$ is a homeomorphism between two sets $X,Y\subset{\mathbb R}^m$, the number of components of $X$ agree with of $Y$. Note that $X$ is  path connected if  there is exactly  one component, that is, $C_p=X$ for all $p\in X$.

Once defined a path component, it would be nice to have a simple technique  that provides a method for computing the number of path components.  One could think that if $X=A\cup B$, with $A\cap B=\emptyset$ and both $A$ and $B$ are path connected, then $A$ and $B$ are, indeed, the components of $X$. This does not hold in general. For example,  take the Euclidean line ${\mathbb R}$. Then ${\mathbb R}=(-\infty,0]\cup (0,\infty)$ is a partition by path connected subsets which obviously are not  the components of ${\mathbb R}$ because ${\mathbb R}$ is path connected.  This example shows that if we have a such partition, we must impose more conditions that ensure that, indeed, they are the path components. The result that we need is:

\begin{method}  Let $X\subset{\mathbb R}^m$ a set in Euclidean $m$-space. A partition of $X$  by path connected open sets in $X$ is, indeed, the partition of the path connected components.
\end{method}
The new ingredient  in the statement  is that the subsets of the partition are \emph{opens in $X$}. This just fails in the above example  because the interval $(-\infty,0]$ is not an open of ${\mathbb R}$: if $(-\infty,0]$ is open, there would exist a $1$-ball centered at $x=0$ of type $(0-r,0+r)$ and included in $(-\infty,0]$, which it is not possible.

\begin{figure}[hbtp]
\begin{center}
\includegraphics{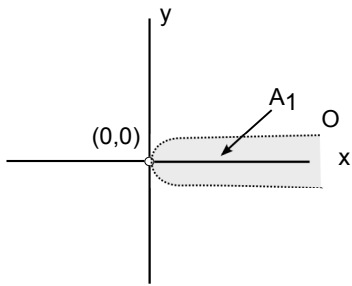}
\end{center}
\caption{The set $X$  formed by the coordinate axis of ${\mathbb R}^2$ except the origin $(0,0)$  has $4$ path components, namely, each one of the  half-axis. Observe that  $A_1=(0,\infty)\times\{0\}$  is an open set in $X$.}\label{fig4}
\end{figure}

We show how to apply our method with the next example. Consider the set formed by the coordinates axis of ${\mathbb R}^2$ except the origin $(0,0)$:
$$X=\big(\left({\mathbb R}\times\{0\}\right)\cup\left(\{0\}\times{\mathbb R}\right)\big)-\{(0,0)\},$$
 see Fig. \ref{fig4}. It is easy to imagine that $X$ has $4$ path components, in fact, the $4$ half-axis. Our method applies as follows. First, each one of the half-axis is path connected   due to it is homeomorphic to an interval, e.g. the interval $(0,\infty)$. The next step is to prove that each half-axis is an open set in $X$. We only do the argument for the half-axis $A_1=(0,\infty)\times\{0\}$. The set $O=(0,\infty)\times{\mathbb R}$ is an open set in ${\mathbb R}^2$ because is the Cartesian product of two opens of ${\mathbb R}$. Then $A_1=X\cap O$ proving that $A_1$ is an open set in $X$. We point out here that not all opens of ${\mathbb R}^2$ are the Cartesian product of opens of ${\mathbb R}$, as for example a $2$-ball, but any Cartesian product of open sets is an open of ${\mathbb R}^2$. The open set $O$ is not unique, as it   appears in Fig. \ref{fig4}, where the set $O$, the shading area, is not a Cartesian product of opens of ${\mathbb R}$.

Once established the method to compute the number of path components, we need the next

\begin{definition} Let $X\subset{\mathbb R}^m$ and $p\in X$. We say that $p$ is an intersection point of order $n\in{\mathbb N}$ if the set $X-\{p\}$ has exactly $n$ path components. This number of components is called the intersection order of $p$ and we denote by  $O(p)$.
\end{definition}

Again the intersection order is a topological property in the sense that if $\phi:X\rightarrow Y$ is a homeomorphism, then $O(p)=O(\phi(p))$ for all $p\in X$.  We compute the intersection orders in explicit examples of sets of ${\mathbb R}^2$.
\begin{enumerate}
\item In the unit circle ${\mathbb S}^1=\{(x,y): x^2+y^2=1\}$ all points have order $1$ because the circle ${\mathbb S}^1$ minus one point is homeomorphic to an interval, which is path connected.
\item All points of a straight-line of ${\mathbb R}^2$ have order $2$. If the straight-line is the axis of abscissas $X=\{(x,0):x\in{\mathbb R}\}$ and $(a,0)\in X$, then
    $$X-\{(a,0)\}= \{(x,0):-\infty<x<a\}\cup \{(x,0):a<x<\infty\}.$$
    Each one of the sets $A=\{(x,0):-\infty<x<a\}$ and $B=\{(x,0):a<x<\infty\}$ are homeomorphic to an interval, namely, $(-\infty,a)$ and $(a,\infty)$, respectively. This means that $A$ and $B$ are path connected. Moreover,   it is clear that $A=X\cap\left((-\infty,a)\times\mathbb R\right)$ and $B=X\cap\left((a,\infty)\times\mathbb R\right)$, being $(-\infty,a)\times\mathbb R$ and   $(a,\infty)\times\mathbb R$   opens of ${\mathbb R}^2$.
\item For the coordinates axis  $X=\left({\mathbb R}\times\{0\}\right)\cup\left(\{0\}\times{\mathbb R}\right)$, we have
    $$O(p)=\left\{\begin{array}{ll} 4& p=(0,0)\\ 2 &p\not=(0,0).\end{array}\right.$$
The case $p=(0,0)$ has been previously treated.  Suppose now that $p\not=(0,0)$. If we suppose that $p=(a,0)$ with $a>0$, then $X-\{p\}=A_1\cup A_2$, where $A_1=\{(x,0):a<x<\infty\}$ and $A_2=X-(\{p\}\cup A_1)$. The set $A_1$ is path connected because it is homeomorphic to the interval $(a,\infty)$. The set $A_2$ is path connected too since is the union of two path connected sets, namely, the $y$-axis (a straight-line) and $\{(x,0):a<x<\infty\}$ (homeomorphic to $(a,\infty)$) and both subsets have non-empty intersection. Finally,   $A_1$ and $A_2$ are open sets in $X-\{p\}$ because
$$A_1=(X-\{p\})\cap \left((a,\infty)\times{\mathbb R}\right),\ A_2=(X-\{p\})\cap \left((-\infty,a)\times{\mathbb R}\right),$$
where $(a,\infty)\times {\mathbb R}$ and $(-\infty,a)\times{\mathbb R}$ are two opens of ${\mathbb R}^2$.
\end{enumerate}

\section{Classifying the letters of the alphabet.}

Let us go back with  the letters of the alphabet of Fig. \ref{fig2}. The work to prove or disprove that  two sets $X,Y\subset {\mathbb R}^2$ are homeomorphic has a different flavor. If we hope that the answer is `yes', then we have to find an \emph{explicit} homeomorphism between both. Following this idea,
by a process of stretching, projecting and rotating, we proved in Fig. \ref{fig1} that {\sf C}, {\sf I} and {\sf L} are homeomorphic.  Working with the allowed topological operations, we find immediately the next five groups (with more than one letter) of homeomorphic letters listed in Fig. \ref{fig5}.
\begin{figure}[hbtp]
\begin{center}
{\sf A R}\\
{\sf C I J L M N S U V W Z}\\
{\sf D O}\\
{\sf E F G T Y}\\
{\sf H K}
\end{center}
\caption{A first step in the topological classification of the letters of the alphabet by using homeomorphisms between pairs of letters.}\label{fig5}
\end{figure}

We point out that in this step of discussion, we can not topologically distinguish, for example, between the letter {\sf C} and the letter {\sf E}. Although our intuition may lead us that  it is not possible to transform {\sf C} to {\sf E} by the admissible deformations (stretching,  contracting and so on), this  only indicates that {\it we can do it explicitly}, but this does not prove that {\sf C} and {\sf E} are not homeomorphic.

\begin{figure}[hbtp]
\begin{center}\includegraphics{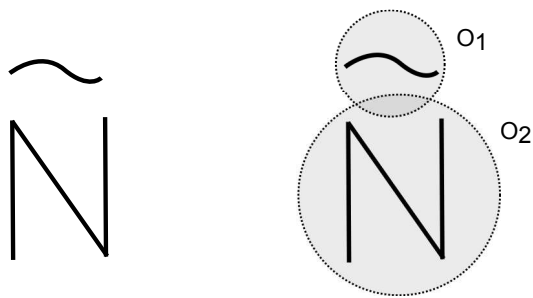}
\end{center}
\caption{The letter {\sf {Ñ}} has two path components. In particular, the letter {\sf {Ñ}} is not homeomorphic to another  letter that appears in Fig. \ref{fig2}, which are all path connected.}\label{fig6}
\end{figure}

By contrast, if we believe that two letters are not homeomorphic,  then we need to give a topological property that satisfies $X$ but not $Y$. Our strategy was announced at the beginning by using the notion of path connectedness and the intersection order. First, all letters of Fig. \ref{fig2} are path connected,  and thus the topological property `to be path connected' does not distinguish  any two letters. Here we may stop a moment and consider the Spanish letter {\sf { Ñ}} (pronounced as \verb|énye|) which does not appear in Fig. \ref{fig2}. This letter has two path components, namely, the symbols {\sf N} and $\sim$. Each one is path connected, indeed, both are homeomorphic to a closed interval $[0,1]$ and both are opens  in the letter {\sf { Ñ}} such as it is depicted  in Fig. \ref{fig6}.

We return with the alphabet of Fig. \ref{fig2}. We know that if  $X\cong Y$ by a homeomorphism $\phi$ and if $p\in X$, then  $O(p)=O(\phi(p))$. But this holds {\it for all} points of $X$. Then we introduce the next notation. Let $n\in{\mathbb N}$ and denote
$$N(n,X)=\mbox{ card}\{p\in X: O(p)=n\}.$$
 Thus, if $X$ is homeomorphic to $Y$, $N(n,X)=N(n,Y)$ for all $n\in{\mathbb N}$. Our work will consist to compute $N(n,X)$ for each letter $X$ of the alphabet. Here we make use of our method. For the case of a letter $X$, we will have
 $N(n,X)=0$ for almost $n\in{\mathbb N}$, and there will also exist $n\in{\mathbb N}$ with $N(n,X)=\infty$. However, there will be some  numbers $n\in{\mathbb N}$ such that  $0<N(n,X)<\infty$.

 \begin{remark} For those letters that have in the typography `segments', the end points of each segment have order $1$. This occurs, for example, with the two lowest points of the letter {\sf A}, or with the two ends points of the letter {\sf I}.
\end{remark}

\begin{figure}[hbtp]
 \begin{center}\includegraphics{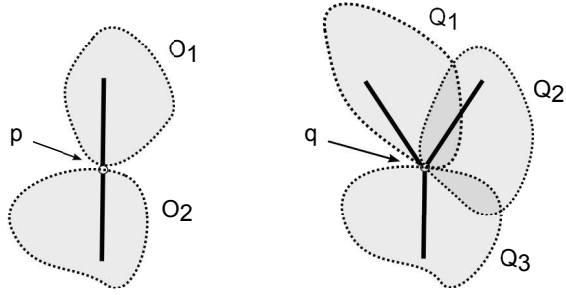}\end{center}
\caption{All points of the letter {\sf I} have order $2$, except the end points. In the letter {\sf Y}, the point $q$ has order $O(q)=3$.}\label{fig7}
 \end{figure}

 First we begin by distinguishing two letters that appear in Fig. \ref{fig5}.  We illustrate the method for the letters {\sf I} and {\sf Y} which are representative of the general case and a similar process may done for any pair of letters.  With the same argument that we did with a straight-line in the previous section, it is clear that $O(p)=2$ for all $p\in{\sf I}$ except the end points which have intersection order $1$. See Fig. \ref{fig7}. Here we use strongly the assumption that the letter is a curve, that is, a one dimensional set, in such a way that deleting a point, we break the letter {\sf I} into two pieces. On the contrary, in the letter {\sf Y} there is one point $q\in{\sf Y}$ such that $O(q)=3$: the point $q\in{\sf Y}$ is just the place where the three segments ensemble to construct the letter {\sf Y}.  In Fig. \ref{fig7},  each component of ${\sf I}-\{p\}$ and ${\sf Y}-\{q\}$ is covered by an open set of ${\mathbb R}^2$ ($O_1$ and $O_2$ for the letter {\sf I} and $Q_1, Q_2,Q_3$ for {\sf Y}).

It also occurs that $N(1,{\sf I})=2$ and $N(1,{\sf Y})=3$, which coincide with the end points of the segments. Finally, we have
 $$N(n,{\sf I})=\left\{ \begin{array}{ll}
 2&n=1\\
 \infty& n=2\\
 0 &n\not=1,2\end{array}\right.,\ \
  N(n,{\sf Y})=\left\{\begin{array}{ll}
  3&n=1\\
  \infty &n=2\\
  1& n=3\\
 0 &n\not=1,2,3.\end{array}\right.$$

This proves definitively that the letters {\sf I} and {\sf Y} are not homeomorphic.

We proceed to compute the intersection order of the letters {\sf O,  A, B, H, P, Q} and {\sf X}. Because the letter {\sf O} is homeomorphic to the unit circle ${\mathbb S}^1$, we conclude from the previous section that
$$N(n,{\sf O})=\left\{\begin{array}{ll}
  \infty &n=1\\
 0 &n\not=1.\end{array}\right.$$
Furthermore, we have $N(n,{\sf B})=N(n,{\sf O})$ for all $n\in{\mathbb N}$  and

$$N(n,{\sf A})=N(n,{\sf P})=\left\{\begin{array}{ll}
 \infty &n=1,2\\
 0 &n\not=1,2\end{array}\right.$$
 $$ N(n,{\sf H})=\left\{\begin{array}{ll}
4 &n=1\\
 \infty &n=2\\
 2 &n=3\\
 0&n\not= 1,2,3\end{array}\right.,\ \ N(n,{\sf Q})=\left\{\begin{array}{ll}
 \infty &n=1,2\\
 1&n=3\\
 0 &n\not=1,2,3\end{array}\right.$$
 $$N(n,{\sf X})=\left\{\begin{array}{ll}
 4 &n=1\\
 \infty &n=2\\
 1&n=4\\
 0 &n\not=1,2,4.\end{array}\right.$$

After these computations,   we obtain in Fig. \ref{fig8} new groups of letters  that, together the ones of Fig. \ref{fig5}, satisfy that two letters of different groups are not homeomorphic.
\begin{figure}[hbtp]
\begin{center}
{\sf A R P}\\
 {\sf  D O B}\\
 {\sf Q}\\
 {\sf X}
 \end{center}
 \caption{A letter in each of the above four groups is not homeomorphic to other letter in other group or in the groups that appear in Fig. \ref{fig5}.}\label{fig8}
 \end{figure}
 Looking these letters it seems intuitive  that the letters {\sf A} and {\sf P} are not homeomorphic, although the intersection orders coincide. The same occurs between {\sf B} and {\sf O}.  In both cases, we need to refine the above arguments. First let us work with {\sf B} and {\sf O}. Assume that there is a homeomorphism $\phi:{\sf B}\rightarrow{\sf O}$ and take a point $p\in {\sf B}$. Then ${\sf B}-\{p\}\cong{\sf O}-\{\phi(p)\}$. In particular, $N(n,{\sf B}-\{p\})=N(n,{\sf O}-\{\phi(p)\})$. The key point $p$ in the letter {\sf B} is indicated in Fig. \ref{fig9}. Let $\phi(p)$ be the corresponding point in the letter {\sf O}. Here all points of {\sf O} have intersection orders are $1$. In the set ${\sf B}-\{p\}$ we compute the intersection order of the point $q$, obtaining $O(q)=3$.  However, in the set ${\sf O}-\{\phi(p)\}$ all points have intersection order $2$, obtaining a contradiction. In other words, and summarizing, if we remove two points $p,q\in {\sf B}$ as in  Fig. \ref{fig8}, then ${\sf B}-\{p,q\}$ should be homeomorphic to ${\sf O}-\{\phi(p),\phi(q)\}$, but the first set has three path components and the second one has only two path components. This proves that  {\sf B} and {\sf O} are not homeomorphic.

\begin{figure}[hbtp]
\begin{center}\includegraphics{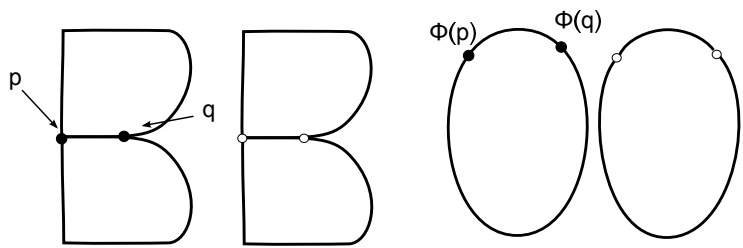}\end{center}
\caption{We remove the points $p$ and $q$ from the letter {\sf B} obtaining $3$ path components. However, for any two points that we delete in the letter {\sf O}, there are always $2$ path components.}\label{fig9}
\end{figure}

\begin{figure}[hbtp]
\begin{center}\includegraphics{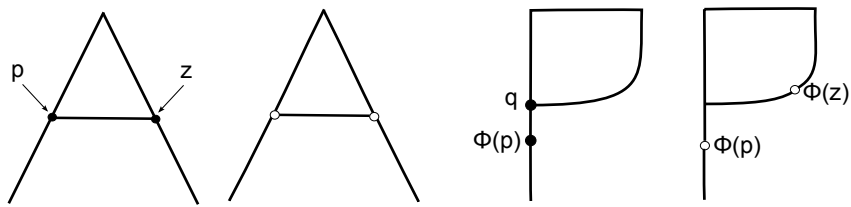}\end{center}
\caption{If we take off the points $p$ and $z$ from the letter {\sf A}, we obtain $4$ components. In the letter {\sf P}, and independently if the point $\phi(p)$ is $q$, if we delete another point, there are $2$ or $3$ components.}\label{fig10}
\end{figure}

Finally, we distinguish the letter {\sf A} from {\sf P}. The argument is similar as above but with a bit difference. Assume that $\phi:{\sf A}\rightarrow{\sf P}$ is a homeomorphism and let $p\in {\sf A}$ be the point indicated in Fig. \ref{fig10}. We know that $O(p)=2$. Then we look for in the letter {\sf P} the points with intersection order $2$. These points belong to the vertical segment below the point $q\in {\sf P}$ that appears in Fig. \ref{fig10}, even possibly the very point $q$. Then ${\sf A}-\{p\}\cong {\sf P}-\{\phi(p)\}$. Now we delete the point $z$ in ${\sf A}-\{p\}$ such it is marked in Fig. \ref{fig10}. Then it remains $4$ path components, that is, $O(z)=4$ in ${\sf A}-\{p\}$. However, if we remove the point $\phi(z)$, the number of path components that remains is  $2$ or $3$, depending on $\phi(p)=q$ or $\phi(p)\not=q$, respectively.  This contradiction proves definitively that {\sf A} is not homeomorphic to {\sf P}.

As a conclusion,

\begin{classification} The topological classification of the letters of the alphabet written in the Sans Serif font of \TeX \ is the following:
\begin{center}
{\sf A R}\\
{\sf B}\\
{\sf C I J L M N S U V W Z}\\
{\sf D O}\\
{\sf E F G T Y}\\
{\sf H K}\\
{\sf P}\\
{\sf Q}\\
{\sf X}
\end{center}
\end{classification}

\section{Going further: other sets and dimensions.}

Having obtained the classification of the letters of the alphabet, one may extend  the same technique in other settings. May we use the concept of intersection order to classify topologically other subsets of ${\mathbb R}^2$? Is it possible to generalize to other dimensions, as for example, in $3$-dimensional Euclidean space ${\mathbb R}^3$?

\begin{figure}[hbtp]
\begin{center}\includegraphics{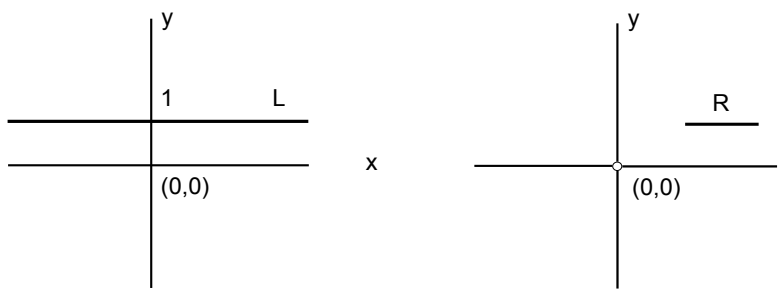}\end{center}
\caption{On the left, the horizontal line $L$ splits  ${\mathbb R}^2$ into two path components but $L$ is homeomorphic to $R=\{(x,1):1<x<2\}$ which, when it is removed from $Y={\mathbb R}^2-\{(0,0)\}$, it remains only $1$ path component.}\label{fig11}
\end{figure}

The idea that lies behind our method is that a letter is a $1$-dimensional set and we are removing $0$-dimensional objects, indeed, points. The problem  in ${\mathbb R}^2$ is more complicated if we take  non $1$-dimensional subsets. For example, let $X={\mathbb R}^2$, the entire plane, and $Y$ the punctured plane $Y={\mathbb R}^2-\{(0,0)\}$.  Both sets are not homeomorphic, but the method of intersection orders seems useless because for any $p\in X$ and $q\in Y$, $O(p)=O(q)=1$.  Even  if we follow with this procedure of deleting point by point,  the  remaining set is path connected. In fact, in ${\mathbb R}^2$ (or ${\mathbb R}^2-\{(0,0)\}$), the complement of an infinite countable set is path connected.   What should we try? As $X$ and $Y$ are $2$-dimensional objects, one may remove $1$-dimensional sets, such as curves.  The simplest case of curve is a straight-line. A possible argument could be the following. Assume that $\phi:X\rightarrow Y$ is a homeomorphism. Consider $L\subset X$ a straight-line, as for example,  $L=\{(x,1):x\in{\mathbb R}\}$ and we remove $L$ from $X$, obtaining $2$ components. Then $\phi(L)$ is other set homeomorphic to a straight-line.
But, what is the shape of $\phi(L)$ as a subset of ${\mathbb R}^2-\{(0,0)\}$? One may think that $\phi(L)$ divides ${\mathbb R}^2-\{(0,0)\}$ into $2$ path components and this occurs if, for example, $\phi(L)$ is other straight-line. In contrast, $\phi(L)$ may be `small' as, for example, any segment of Euclidean plane. So, the subset $\{(x,1):1<x<2\}$ is homeomorphic to $L$, obtaining that $Y-\phi(L)$ has only $1$ component! See Fig. \ref{fig11}.

\begin{figure}[hbtp]
\begin{center}\includegraphics{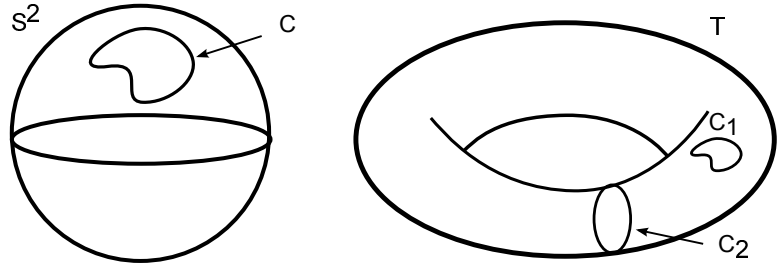}\end{center}
\caption{In sphere ${\mathbb S}^2$ the  curve $C$ divides ${\mathbb S}^2$ in  $2$ path components. On  the doughnut ${\mathbb T}$ there are simple closed curves, as $C_1$, that splits ${\mathbb T}$ into  $2$ components but, in contrast,  if we remove the meridian $C_2$, we obtain a path connected set.}\label{fig12}
\end{figure}

In $3$-dimensional Euclidean space ${\mathbb R}^3$ we may study if  a round sphere ${\mathbb S}^2$ is homeomorphic to the surface of a  doughnut ${\mathbb T}$, called a torus (in topological language we say `is homeomorphic to' ${\mathbb S}^1\times{\mathbb S}^1$). In ${\mathbb S}^2$ there are closed curves $C$ without self-intersections  (or simple closed curves) that when it is deleted from ${\mathbb S}^2$, it rests $2$ path components. In contrast, in ${\mathbb T}$ there are curves, as $C_1$ in Fig. \ref{fig12} with the same property, but there is also curves the equator or a meridian whose complement is path connected. Even more,  we must point out to the reader that  the statement \emph{the complement of a simple closed curve in a sphere has exactly two components}, is true but not obvious! This is the famous Jordan curve theorem, which it is a deep result in topology and its the proof is far to be `trivial' \cite{pe,th}. In $m$-dimensional Euclidean space the topological properties needed to distinguish  objects such as ${\mathbb S}^2$ from ${\mathbb T}$ use more sophisticated properties than path connectedness, but this problem carries out of the initial plan of this article.

\begin{thebibliography}{1}
\bibitem{ch} Gustave Choquet,   \textit{Topology}, Academic Press,  New York, 1966.
\bibitem{mu} James R. Munkres, \textit{Topology}, Prentice-Hall, 2nd ed. Upper Saddle River, 2000.

\bibitem{pe} R. N. Pederson, The Jordan curve theorem for piecewise smooth curves,
\emph{Amer. Math. Monthly} \textbf{76} (1969) 605--610.
\bibitem{th} Carsten Thomassen, The Jordan-Sch\"{o}nflies theorem and the classification of surfaces,
\emph{Amer. Math. Monthly} \textbf{99} (1992) 116--130.
 \bibitem{wiki} Topology. (2014, June 25). In Wikipedia, The Free Encyclopedia, from http://en.wikipedia.org/w/index.php?title=Topology\&oldid=614319857, [Online; accessed 27-June-2014].





 \end{thebibliography}
\end{document}